\documentclass[12pt]{amsart}

\usepackage{amsmath,amssymb,amsthm,color}

\usepackage[margin=1in]{geometry}

\title{Unexpected biases in the distribution of consecutive primes}

\author{Robert J. Lemke Oliver}
\address{Department of Mathematics, Stanford University}
\email{rjlo@stanford.edu}
\address{Department of Mathematics, Tufts University}
\email{robert.lemke\_oliver@tufts.edu}

\author{Kannan Soundararajan}
\address{Department of Mathematics, Stanford University}
\email{ksound@stanford.edu}

\dedicatory{To Professor R. Balasubramanian on his sixty fifth birthday} 

\numberwithin{equation}{section}
\newtheorem{theorem}{Theorem}

\newtheorem{conjecture}[theorem]{Conjecture}
\newtheorem*{mainconjecture}{Main Conjecture}
\newtheorem{proposition}[theorem]{Proposition}
\numberwithin{theorem}{section}
\theoremstyle{remark} 
\theoremstyle{remark} 
\theoremstyle{remark} \newtheorem*{example}{Example}

\renewcommand{\pmod}[1]{\,\left(\text{mod }#1\right)}

\begin{document}

\begin{abstract}  While the sequence of primes is very well distributed in the reduced residue classes $\pmod q$, the distribution of pairs of consecutive primes among the permissible $\phi(q)^2$ pairs of reduced residue classes $\pmod q$ is surprisingly erratic.   This paper proposes a conjectural explanation for this phenomenon, based on the Hardy-Littlewood conjectures.  The conjectures are then compared to numerical data, and the observed fit is very good. 
\end{abstract}

\maketitle

\section{Introduction}

\noindent The prime number theorem in arithmetic progressions shows that the sequence of primes is equidistributed among the reduced residue classes $\pmod q$.   
If the Generalized Riemann Hypothesis is true, then this holds in the more precise form 
\[
\pi(x;q,a) = \frac{\mathrm{li}(x)}{\phi(q)} + O(x^{1/2+\epsilon}), \quad \mbox{where } \mathrm{li}(x) := \int_2^x \frac{dt}{\log t}, 
\]
and $\pi(x;q,a)$ denotes the number of primes up to $x$ lying in the reduced residue class $a\pmod q$.    
Nevertheless it was noticed by Chebyshev that certain residue classes 
seem to be slightly preferred: for example, among the first million primes, we find that
\[
\pi(x_0;3,1) = 499{,}829 \quad \mbox{ and } \quad \pi(x_0;3,2) = 500{,}170, \quad \pi(x_0) = 10^6.  
\]
Chebyshev's bias is beautifully explained by the work of Rubinstein and Sarnak \cite{RubinsteinSarnak} (see \cite{GranvilleMartin} for a survey of related work) who showed (in a certain sense and under some natural conjectures) that $\pi(x;3,2)>\pi(x;3,1)$ for $99.9\%$ of all positive $x$.

What happens if we consider the patterns of residues $\pmod q$ among strings of consecutive primes?   Let $p_n$ denote the sequence of primes in ascending order.  
Let $r \ge 1$  be an integer, and let $\mathbf{a} = (a_1, a_2, \dots, a_{r})$ denote an $r$-tuple of reduced residue classes $\pmod q$.   Define 
\[
\pi(x;q,\mathbf{a}) := \#\{p_n \le x : p_{n+i-1}\equiv a_i \pmod{q} \text{ for each $1\leq i \leq r$}\}, 
\]
which counts the number of occurrences of  the pattern $\mathbf{a} \pmod{q}$ among $r$ consecutive primes the least of which is below $x$.  When $r\ge 2$, little is known about the distribution of such patterns among the primes.  When $r=2$ and $\phi(q) =2$ (thus $q=3$, $4$, or $6$), Knapowski and Tur{\' a}n \cite{KnapowskiTuran} observed that 
all the four possible patterns of length $2$ appear infinitely many times.   The main significant result in this direction is due to D. Shiu \cite{Shiu} who established that for any $q\ge 3$, a reduced residue class $a\pmod q$, and any $r\ge 2$, the pattern $(a, a, \ldots, a)$ occurs  infinitely often.   Recent progress in sieve theory has led to a new 
proof of Shiu's result (see \cite{BFTB}), and moreover in this particular situation Maynard \cite{Maynard} has  shown that $\pi(x;q,(a,\ldots,a)) \gg \pi(x)$.  
 
 Despite the lack of understanding of $\pi(x;q, \mathbf{a})$, any model based on the randomness of the primes would suggest strongly that every permissible 
 pattern of $r$ consecutive primes appears roughly equally often:  that is, if ${\mathbf a}$ is an $r$-tuple of reduced residue classes $\pmod q$, then $\pi(x;q,{\mathbf a}) \sim 
  \pi(x)/\phi(q)^r$.   However, a look at the data might shake that belief!  For example, among the first million primes (for convenience restricting to those greater than 3) we find
\begin{align*}
\pi(x_0;3,(1,1)) = 215{,}873, \quad & \pi(x_0;3,(1,2)) = 283{,}957, \\ 
\pi(x_0;3,(2,1)) = 283{,}957, \quad & \pi(x_0;3,(2,2)) = 216{,}213.
\end{align*}
These numbers show substantial deviations from the expectation that all four quantities should be roughly $250{,}000$.  Further, Chebyshev's bias $\pmod 3$ might have suggested a slight preference for the pattern $(2,2)$ over the other possibilities, and this is clearly not the case.    

The discrepancy observed above persists for larger $x$, and also exists for other moduli $q$.  For example, among the first hundred million primes modulo $10$, there is substantial deviation from the prediction that each of the 16 pairs $(a,b)$ should have about $6.25$ million occurrences.  Specifically, with $\pi(x_0)=10^8$, we find the following.

\begin{center}
\begin{tabular}{ll|r}
$a$ & $b$ & $\pi(x_0;10,(a,b))$ \\
\hline 
1 & 1 & 4{,}623{,}042 \\
  & 3 & 7{,}429{,}438 \\
  & 7 & 7{,}504{,}612 \\
  & 9 & 5{,}442{,}345 \\
3 & 1 & 6{,}010{,}982 \\
  & 3 & 4{,}442{,}562 \\
  & 7 & 7{,}043{,}695 \\
  & 9 & 7{,}502{,}896
\end{tabular} \quad\quad
\begin{tabular}{ll|r}
$a$ & $b$ & $\pi(x_0;10,(a,b))$ \\
\hline 
7 & 1 & 6{,}373{,}981 \\
  & 3 & 6{,}755{,}195 \\
  & 7 & 4{,}439{,}355 \\
  & 9 & 7{,}431{,}870 \\
9 & 1 & 7{,}991{,}431 \\
  & 3 & 6{,}372{,}941 \\
  & 7 & 6{,}012{,}739 \\
  & 9 & 4{,}622{,}916
\end{tabular}
\end{center}
Apart from the fact that the entries vary dramatically (much more than in Chebyshev's bias), the key feature to be observed in this data is that the diagonal classes $(a,a)$ occur significantly less often than the non-diagonal classes.    Chebyshev's bias $\pmod{10}$ states that the residue classes $3$ and $7 \pmod {10}$ very often contain slightly more primes than the residue classes $1$ and $9 \pmod{10}$, but curiously in our data the patterns $(3,3)$ and $(7,7)$ appear less frequently than $(1,1)$ and $(9,9)$; this suggests again that a different phenomenon is at play here.  

The purpose of this paper is to develop a heuristic, based on the Hardy-Littlewood prime $k$-tuples conjecture, which explains the biases seen above.  We are led to conjecture that while the primes counted by $\pi(x;q,\mathbf{a})$ do have density $1/\phi(q)^r$ in the limit, there are large secondary terms in the asymptotic formula which create biases toward and against certain patterns.   The dominant factor in this bias is determined by the number of $i$ for which $a_{i+1}\equiv a_i \pmod{q}$, but there are also lower order terms that do not have an easy description. 

\begin{mainconjecture}
\label{conj:main}
With notation as above, we have
\[
\pi(x;q,\mathbf{a}) = \frac{\mathrm{li}(x)}{\phi(q)^r} \Big(1 + c_1(q;\mathbf{a}) \frac{\log\log x}{\log x} + c_2(q;\mathbf{a})\frac{1}{\log x} + O\Big(\frac{1}{(\log x)^{7/4}}\Big)\Big),
\]
where
\[
c_1(q;\mathbf{a}) = \frac{\phi(q)}{2} \Big(\frac{r-1}{\phi(q)}-\#\{1\leq i < r : a_i \equiv a_{i+1}\pmod{q}\}\Big),
\]
and when $r=2$ the constant $c_2(q;\mathbf{a})$ is given in \eqref{eqn:c2_final}, while if $r\ge 3$
\[
c_2(q;\mathbf{a}) = \sum_{i=1}^{r-1} c_2(q;(a_i,a_{i+1})) + \frac{\phi(q)}{2}\sum_{j=1}^{r-2} \frac{1}{j}\Big(\frac{r-1-j}{\phi(q)} - \#\{i : a_i \equiv a_{i+j+1} \pmod{q}\}\Big).
\]
\end{mainconjecture}

In general, the quantity $c_2(q;\mathbf{a})$ seems complicated, but there are some situations where it simplifies.  
For example, if $\mathbf{a}=(a,a)$ for a reduced residue class $a\pmod q$, then regardless of the choice of $a$ we have  
\begin{equation}
\label{eqn:diagonal_constant}
c_2(q;(a,a))=\frac{\phi(q)\log(q/2\pi) + \log 2\pi}{2} - \frac{\phi(q)}{2} \sum_{p\mid q}\frac{\log p}{p-1}. 
\end{equation}
We can also show that $c_2(q;(a,b)) = c_2(q;(-b,-a))$ for any two reduced residue classes $a$ and $b\pmod{q}$.   Moreover, 
while $c_2(q;(a,b))$ seems involved, the symmetric quantity $c_2(q;(a,b))  + c_2(q;(b,a))$ simplifies nicely: for distinct reduced residue 
classes $a$, $b \pmod{q}$ we have  
\begin{equation} 
\label{1.2} 
c_2(q;(a,b)) + c_2(q;(b,a)) = \log (2\pi )  - \phi(q) \frac{\Lambda(q/(q,b-a))}{\phi(q/(q,b-a))}, 
\end{equation} 
where $\Lambda$ denotes the von Mangoldt function.  In particular, this expression depends only on the difference $b-a$.   

\begin{conjecture}  If $a$ and $b$ are distinct reduced residue classes $\pmod q$, then $\pi(x;q,(a,b)) + \pi (x;q,(b,a))$ equals 
$$ 
2\frac{\mathrm{li}(x)}{\phi(q)^2}  \Big( 1 + \frac{\log \log x}{2\log x} + \Big(\log (2\pi) - \phi(q) \frac{\Lambda(q/(q,b-a))}{\phi(q/(q,b-a))}\Big) \frac{1}{2\log x} + O\Big(\frac{1}{(\log x)^{7/4}} \Big)\Big),  
$$ 
whereas $\pi(x;q,(a,a))$ equals 
$$ 
\frac{\mathrm{li}(x)}{\phi(q)^2} \Big( 1 - \frac{\phi(q)-1}{2} \frac{\log \log x}{\log x}  + \Big(\phi(q)\log \frac{q}{2\pi} +\log 2\pi -\phi(q) \sum_{p|q} \frac{\log p}{p-1}\Big) 
\frac{1}{2\log x} + O\Big(\frac{1}{(\log x)^{7/4}} \Big)  \Big). 
$$
\end{conjecture}

We give a few amusing consequences of the Main Conjecture. 
The famous biases $\pi(x) < \mathrm{li}(x)$, or $\pi(x;3,1) < \pi(x;3,2)$, or $\pi(x;4,1) < \pi(x;4,-1)$ are known to be false infinitely often.  However we conjecture that the robust biases in pairs of consecutive primes $\pmod 3$ or $\pmod 4$ may hold always and from the very start!

\begin{conjecture}
\label{conj:always}
Let $q=3$ or $4$, and let $a$ be either $1\pmod q$ or $-1\pmod q$.   Then for all 
$x\geq 5$,  we have $\pi(x;q,(a,-a)) > \pi(x;q,(a,a))$.  Indeed for large $x$ we have 
$$ 
\pi(x;q,(a,-a)) - \pi(x;q,(a,a)) = \frac{x}{4(\log x)^2} \log \Big( \frac{2\pi }{q}\log x \Big) + O\Big( \frac{x}{(\log x)^{11/4}}\Big).
$$ 
\end{conjecture}

Given a prime $q$, the product of two consecutive primes prefers to be a quadratic non-residue rather than a quadratic residue. 

\begin{conjecture} 
\label{quad} 
Let $q$ be a fixed odd prime.  For large $x$ we have 
$$ 
\sum_{p_n\le x} \Big(\frac{p_n}{q}\Big) \Big(\frac{p_{n+1}}{q}\Big) = - \frac{x}{2(\log x)^2} \log\Big(\frac{2\pi \log x}{q}\Big) 
+ O \Big( \frac{x}{(\log x)^{11/4}}\Big). 
$$ 
\end{conjecture} 

The constants in the Main Conjecture also simplify dramatically if one only cares about patterns exhibited by $p_n$ and $p_{n+k}$ for $k\geq 2$.

\begin{conjecture} \label{conj:skip}
If $k\geq 2$ and $a$ and $b$ are distinct reduced residues $\pmod{q}$, then
\[
\#\{p_n\leq x: p_n \equiv a\pmod{q}, p_{n+k} \equiv b\pmod{q} \} = \frac{\mathrm{li}(x)}{\phi(q)^2}\Big(1 + \frac{1}{2(k-1)}\frac{1}{\log x} + O \Big( \frac{1}{(\log x)^{7/4}}\Big)\Big),
\] 
while
\[
\#\{p_n\leq x: p_n \equiv p_{n+k} \equiv a\pmod{q}\} = \frac{\mathrm{li}(x)}{\phi(q)^2}\Big(1 - \frac{\phi(q)-1}{2(k-1)}\frac{1}{\log x} + O \Big( \frac{1}{(\log x)^{7/4}}\Big)\Big).
\]
\end{conjecture}

Form a $\phi(q)\times \phi(q)$ transition matrix (with rows and columns indexed by reduced residue classes) and the $(a,b)$-th entry being the 
probability that a prime $p_{n}\equiv a\pmod{q}$ is followed by $p_{n+1}\equiv b\pmod{q}$.  
Then Conjecture \ref{conj:skip} shows that the corresponding transition matrix going from $p_n$ to $p_{n+2}$ is not the square of the transition matrix going from $p_n$ 
to $p_{n+1}$.  Thus the primes $\pmod{q}$ are not Markovian, and this may also be seen directly from the Main Conjecture by the formula given for $c_2(q;\mathbf{a})$ 
when $r\geq 3$ (which is used to derive Conjecture \ref{conj:skip}).

The ideas that lead to the Main Conjecture imply that there will be symmetries between the number of occurrences of different patterns.

\begin{conjecture}
\label{conj:symmetry}
Given $\mathbf{a}$ and $q$ as above, define $\mathbf{a}^{\mathrm{opp}} = (-a_r,-a_{r-1},\dots,-a_1)$.  For large $x$ we have
\[
\pi(x;q,\mathbf{a}) = \pi(x;q,\mathbf{a}^\mathrm{opp}) + O(x^{1/2+\epsilon}).
\]
\end{conjecture}

\begin{example}
We find $$\pi(10^{11};7,(1,6,3))=24{,}344{,}117$$ and $$\pi(10^{11};7,(4,1,6))=24{,}349{,}025,$$ while the nearest number of occurrences of another pattern is $$\pi(10^{11};7,(6,2,1))=24{,}570{,}765.$$
\end{example}

If the modulus is a prime power, there are additional symmetries.

\begin{conjecture}
\label{conj:primepower}
Let $q$ be a prime and let $v\geq 2$.  If $\mathbf{a}=(a_1,\dots,a_r)$ and $\mathbf{b}=(b_1,\dots,b_r)$ are such that $a_1 \equiv b_1 \pmod{q}$ and $a_{i+1}-a_i \equiv b_{i+1}-b_i \pmod{q^v}$ for each $1\leq i <r$, then
\[
\pi(x;q^v,\mathbf{a}) = \pi(x;q^v,\mathbf{b}) + O(x^{1/2+\epsilon}).
\]
In particular, if $a$ is odd, then, up to an error $O(x^{1/2+\epsilon})$, $\pi(x;2^v,(a,b))$ depends only on $b-a \pmod{2^v}$.
\end{conjecture}
\begin{example}
We find
\begin{align*}
\pi(10^{11};8,(1,3))=278{,}676{,}326,& & \pi(10^{11};8,(3,5))= 278{,}696{,}997,\\
\pi(10^{11};8,(5,7))= 278{,}692{,}843, &\quad\text{ and}& \pi(10^{11};8,(7,1))=278{,}681{,}776.
\end{align*}
\end{example}

In the direction of these conjectures, the earliest work we found is the paper of Knapowski and Tur{\' a}n \cite{KnapowskiTuran} who ``guess" that the events $p_{n}\equiv a \pmod 4$ and $p_{n+1} \equiv b \pmod 4$ for the four possibilities of $a$ and $b$ are ``not equally probable."  However Knapowski and Tur{\' a}n go on to suggest that $\pi(x;4,(1,1)) = o(\pi (x))$, which is now definitively false by Maynard's work \cite{Maynard}.   The paper \cite{KnapowskiTuran} was published after the death of both authors, and perhaps they had something else in mind, maybe along the lines of our Conjecture \ref{conj:always} above? More recently, in Ko \cite{Ko} numerical results observing the biases in the distribution of consecutive primes for small moduli are given.  The paper by Ash, Beltis, Gross and Sinnott \cite{ABGS} again observes these biases in pairs of consecutive primes and initiates an attempt toward understanding them based on the Hardy-Littlewood conjectures.  The heuristic expression in \cite{ABGS} is a large sum of singular series, and as the authors note, it is unclear from that expression whether $\pi(x;q,(a,b))$ tends to $\pi(x)/\phi(q)^2$ for large $x$.  They also note symmetries akin to Conjectures \ref{conj:symmetry} and \ref{conj:primepower} for pairs of consecutive primes.


In the Main Conjecture we expect that the remainder term $O\big((\log x)^{-7/4}\big)$ is given by a 
sum involving the zeros of Dirichlet $L$-functions $\pmod q$.  The main terms given in the Main Conjecture are the same for all repeating patterns $(a, a, \ldots, a)$; nevertheless numerically one observes some deviations in the counts of such patterns, and we expect the lower order fluctuations to account for these deviations.  In addition to the contributions from zeros, which we expect to be oscillating, there also appear to be non-oscillating lower order terms of size $(\log \log x/\log x)^2$, which may play a bigger role for the computable ranges of $x$.  We hope to understand these lower order terms in future work.     
  
  An initial guess for why there is a bias against the repeating patterns might be that, after a prime occurs that is $a \pmod{q}$, all other classes have a chance to represent a prime before $a$ occurs again.  However, a straightforward application of the Selberg sieve shows that the number of primes for which $p_{n+1}-p_n < q$ is $O(x/\log^2 x)$, which is of a smaller order of magnitude than the bias predicted by the Main Conjecture.

Though we do not pursue this here, it should be possible to prove unconditional analogues of the Main Conjecture in other settings, for example to numbers free of small prime factors or for squarefree integers (in the latter case, the biases will be manifested already at the level of the constant in the main term).  More generally, analogous biases seem to arise for many other sifted sets, for example in the sums of two squares.  We also mention two other settings in which large biases are seen:  the distribution of prime geodesics for compact hyperbolic surfaces into various homology classes (see the discussion at the end of \cite{RubinsteinSarnak}),  and the  recent work of Dummit, Granville, and Kisilevsky \cite{DGK} concerning the distribution of numbers that are products of two primes.

\noindent {\bf Acknowledgements.}  The first author is partially supported by an NSF postdoctoral fellowship, DMS 1303913.  The second author is partially supported by  the NSF, and a Simons Investigator Award from the Simons Foundation.    We would like to thank Tadashi Tokieda whose lecture on ``Rock, paper, scissors in probability" inspired the present work, James Maynard for drawing our attention to \cite{KnapowskiTuran}, Paul Abbott for pointing us to \cite{Ko}, and Alexandra Florea, Andrew Granville and Peter Sarnak for helpful comments.   

\section{The heuristic for $r=2$}
\label{sec:heuristic}

\noindent In this section we develop a heuristic explanation of the Main Conjecture in the case $r=2$.   The heuristic (like several other conjectures about the primes, see for example \cite{Gallagher, GoldstonLedoan, GLR, MontgomerySound, ORW}) is based upon the Hardy-Littlewood prime $k$-tuples conjecture.    We begin by reviewing quickly the Hardy-Littlewood conjectures and some related results, before proceeding to develop an analogue suitable for understanding $\pi(x;q,{\mathbf a})$.  

\noindent {\bf The Hardy-Littlewood conjectures.}   Let $\mathcal{H}$ be a finite subset of $\mathbf{Z}$ and let $\mathbf{1}_\mathcal{P}$ 
   denote the characteristic function of the primes.  In a strong form, the Hardy-Littlewood conjecture asserts that
\[
\sum_{n \le x} \prod_{h\in\mathcal{H}} \mathbf{1}_\mathcal{P}(n+h) = \mathfrak{S}(\mathcal{H}) \int_2^x 
\frac{dy}{(\log y)^{|{\mathcal{H}}|}} + O(x^{1/2+\epsilon}),
\]
where the singular series $\mathfrak{S}(\mathcal{H})$ is given by 
\[
\mathfrak{S}(\mathcal{H}) = \prod_p \Big(1-\frac{\# (\mathcal{H}\, \mathrm{mod}\, p)}{p}\Big)\Big(1-\frac{1}{p}\Big)^{-|\mathcal{H}|}.
\]
In our calculations, it will be important to understand the behavior of the singular series ``on average."   Here 
Gallagher \cite{Gallagher} established that for any $k\ge 1$ and as $h\to \infty$,
\begin{equation}
\label{eqn:gallagher}
\sum_{\begin{subarray}{c} \mathcal{H}\subseteq [1,h] \\ |\mathcal{H}|=k \end{subarray}} \mathfrak{S}(\mathcal{H}) \sim \binom{h}{k} \sim \frac{h^k}{k!}, 
\end{equation}
so that the singular series is $1$ on average.  A refined version of this asymptotic was established by Montgomery and Soundararajan \cite{MontgomerySound}, who introduced the modified singular series 
$$ 
{\mathfrak S}_0({\mathcal H}) = \sum_{{\mathcal T}\subset {\mathcal H}} (-1)^{|{\mathcal H} \setminus \mathcal{T}|}\mathfrak{S}(\mathcal{T}), \qquad  \text{so that} \qquad {\mathfrak S}({\mathcal H}) = \sum_{\mathcal{T} \subset {\mathcal H}} {\mathfrak S}_0({\mathcal T}), 
$$ 
with $ {\mathfrak S}(\emptyset) = {\mathfrak S}_0(\emptyset) = 1$.   The modified singular series ${\mathfrak S}_0$ 
arises naturally in the following version of  the Hardy-Littlewood conjecture (thinking of the elements of ${\mathcal H}$ as being small in comparison to $x$): 
$$ 
\sum_{n\le x} \prod_{h \in \mathcal{H}}\Big ( \mathbf{1}_{\mathcal{P}} (n+h) - \frac{1}{\log n} \Big)  
= {\mathfrak S}_0({\mathcal H})  \int_2^x \frac{dy}{(\log y)^{|{\mathcal H}|}} + O(x^{1/2+\epsilon}), 
$$ 
and the term $1/\log n$ that is subtracted above arises naturally as the probability that the ``random number" $n+h$ is prime.   Montgomery and Soundararajan showed that 
\begin{equation}
\label{eqn:montgomery-sound}
\sum_{\begin{subarray}{c} \mathcal{H}\subseteq [1,h] \\ |\mathcal{H}|=k \end{subarray}} \mathfrak{S}_0(\mathcal{H}) = \frac{\mu_k}{k!} (-h\log h +Ah)^{k/2}+O_k(h^{k/2-1/(7k)+\epsilon}),
\end{equation}
where $\mu_k$ is the $k$-th moment of the standard Gaussian (in particular, $\mu_k=0$ if $k$ is odd) and $A$ is a constant independent of $k$.   This refines Gallagher's asymptotic \eqref{eqn:gallagher}, and shows that ${\mathfrak S}_0({\mathcal H})$ exhibits roughly square-root cancelation in each variable.

\noindent {\bf Modified Hardy-Littlewood conjectures.}  We need a slight modification of the Hardy-Littlewood conjecture, taking into account congruence conditions $\pmod q$.  For any integer $q\geq 1$ and a finite subset $\mathcal{H}$ of the integers, we define the singular series at the primes away from $q$ by 
\[
\mathfrak{S}_q(\mathcal{H}) := \prod_{p\nmid q} \Big(1 - \frac{\#(\mathcal{H}\,\mathrm{mod}\, p) }{p}\Big)\Big(1-\frac{1}{p}\Big)^{-|\mathcal{H}|}. 
\]
 If $a\pmod{q}$ is such that $(h+a,q)=1$ for all $h\in\mathcal{H}$,  then we expect that
\begin{equation}
\label{eqn:hardy-littlewood}
\sum_{\begin{subarray}{c} n<x \\ n \equiv a\pmod{q} \end{subarray}} \prod_{h\in\mathcal{H}} \mathbf{1}_{\mathcal{P}}(n+h) \sim \mathfrak{S}_q(\mathcal{H}) \Big(\frac{q}{\phi(q)}\Big)^{|\mathcal{H}|} \ \frac{1}{q} \int_2^x \frac{dy}{(\log y)^{|\mathcal{H}|} },
\end{equation}
where the factor $(q/\phi(q))^{|\mathcal{H}|}$ arises because $h+a$ is conditioned to be coprime to $q$ for all $h\in {\mathcal H}$, and the factor $1/q$ arises since we are restricting $n$ to one residue class $\pmod q$.  In analogy with ${\mathfrak S}_0$, it is also useful to define $
\mathfrak{S}_{q,0}(\mathcal{H}) := \sum_{\mathcal{T}\subseteq \mathcal{H}} (-1)^{|\mathcal{H}\setminus \mathcal{T}|} \mathfrak{S}_q(\mathcal{T}),
$
so that
$\mathfrak{S}_q(\mathcal{H}) = \sum_{\mathcal{T}\subseteq \mathcal{H}} \mathfrak{S}_{q,0}(\mathcal{T})$. 
Once again the quantity ${\mathfrak S}_{q,0}$ arises naturally in the asymptotic (conditioning $(h+a,q)=1$ for all $h\in {\mathcal H}$) 
\begin{equation} 
\label{2.4} 
\sum_{\substack{n\le x\\ n\equiv a\pmod{q}}} \prod_{h \in {\mathcal H}} \Big( \mathbf{1}_{\mathcal P}(n+h) - \frac{q}{\phi(q) \log n }\Big) \sim  \mathfrak{S}_{q,0}(\mathcal{H}) \Big(\frac{q}{\phi(q)}\Big)^{|\mathcal{H}|} \ \frac{1}{q} \int_2^x \frac{dy}{(\log y)^{|\mathcal{H}|} },  
\end{equation} 
 where the term $q/(\phi(q)\log n)$ being subtracted arises naturally as the probability that $n+h$ is prime, conditioned on the fact that $n+h$ is coprime to $q$.  
 
\noindent {\bf First steps toward the conjecture.}  Let $a$ and $b$ be two reduced residue classes $\pmod q$, and let $h$ be a positive integer with $h\equiv b-a \pmod q$.  We now formulate a conjecture for the number of primes $n \le x$ with $n\equiv a\pmod q$ and such 
that the next prime after $n$ is $n+h$.  The gaps between consecutive primes are conjectured to be 
distributed like a Poisson process with mean $\sim \log x$ (and Gallagher showed that this follows from the Hardy-Littlewood conjectures), and so $h$ should be thought of as a parameter on 
the scale of $\log x$.  With this in mind, we are interested in 
\begin{align}
\label{2.5}
&\sum_{\substack{ n\le x \\ n\equiv a\pmod q}} \mathbf{1}_{\mathcal P}(n)  \mathbf{1}_{\mathcal P}(n+h) 
\prod_{\substack{0<t< h \\ (t+a,q)=1}} \Big(1 -  \mathbf{1}_{\mathcal P}(n+t) \Big) \nonumber \\
&= 
\sum_{\substack{ n\le x \\ n\equiv a\pmod q}} {\mathbf 1}_{\mathcal P}(n) \mathbf{1}_{\mathcal P}(n+h)\prod_{\substack{0<t< h \\ (t+a,q)=1}} 
\Big( 1- \frac{q}{\phi(q) \log (n+t)} - \tilde{\mathbf{1}}_{\mathcal P}(n+t) \Big), 
\end{align} 
where, for a variable $n$ conditioned to be coprime to $q$, we set $\tilde{\mathbf{1}}_{\mathcal P}(n) = \mathbf{1}_{\mathcal P} (n) - q/(\phi(q)\log n)$.   Write also $\mathbf{1}_{\mathcal P}(n) = q/(\phi(q)\log n) + \tilde{\mathbf{1}}_{\mathcal P}(n)$ and similarly for ${\mathbf 1}_{\mathcal P}(n+h)$, and then expand out the product in \eqref{2.5}: thus we arrive at (ignoring the small differences between $\log n$, $\log (n+h)$ or $\log (n+t)$)
\begin{equation} 
\label{2.6} 
\sum_{{\mathcal A} \subset \{0, h\}} \sum_{\substack{{\mathcal T} \subset [1, h-1]  \\ (t+a,q)=1 \forall t\in{\mathcal T}}} (-1)^{|{\mathcal T}|} 
\sum_{\substack{ n\le x\\ n\equiv a\pmod q}} \Big(\frac{q}{\phi(q) \log n} \Big)^{2-|{\mathcal A}|}
 \prod_{\substack{t\in [1,h-1] \\ (t+a,q)=1 \\ t\notin {\mathcal T}} } \Big(1-\frac{q}{\phi(q)\log n} \Big) 
 \prod_{t\in {\mathcal A} \cup {\mathcal T}} \tilde{\mathbf 1}_{\mathcal P} (n+t) . 
\end{equation} 
Given reduced residue classes $a$ and $b$, and a positive $h\equiv b-a \pmod q$, we may write 
\begin{equation} 
\label{2.7} 
\#\{ 0< t<h: \ (t+a,q)=1 \} = \frac{\phi(q)}{q} h + \epsilon_q(a,b), 
\end{equation} 
where $\epsilon_q(a,b)$ is independent of $h$.  We also write for convenience 
\begin{equation} 
\label{2.8} 
\alpha(y) = 1 -\frac{q}{\phi(q)\log y}. 
\end{equation} 
Appealing now to the conjectured relation \eqref{2.4}, we are led to hypothesize that the quantity in \eqref{2.5} (and \eqref{2.6}) is 
\begin{equation} 
\label{2.9}
\sim \sum_{{\mathcal A} \subset \{0, h\}} \sum_{\substack{{\mathcal T} \subset [1, h-1]  \\ (t+a,q)=1 \forall t\in{\mathcal T}}} (-1)^{|{\mathcal T}|} {\mathfrak S}_{q,0}({\mathcal A}\cup {\mathcal T}) \Big( \frac{1}{q} \int_2^x \Big(\frac{q}{\phi(q) \log y}\Big)^{2+|{\mathcal T}|} \alpha(y)^{h\phi(q)/q + \epsilon_q(a,b)-|{\mathcal T}|}  dy\Big) . 
\end{equation}

 Before proceeding further, a few points are in order.  Note that $\alpha(x)^{h\phi(q)/q}$ is about $e^{-h/\log x}$, 
 and this exponential decay in $h$ is in keeping with the conjecture that gaps between consecutive primes are 
 distributed like a Poisson process.  Secondly, by replacing ${\mathcal A}$ and ${\mathcal T}$ above with $h-{\mathcal A}$ and $h-{\mathcal T}$, and noting also that $\epsilon_q(a,b)= \epsilon_q(-b,-a)$ we may see that the 
 quantity \eqref{2.9} above does not change if we replace $(a,b)$ by $(-b,-a)$; this is an example of the symmetry between $\pi(x;q, {\mathbf a})$ and $\pi(x;q,{\mathbf a}^{\text{opp}})$ noted in Conjecture \ref{conj:symmetry}.  Similarly, under the hypotheses of Conjecture \ref{conj:primepower}, the conditions satisfied by $h$ and $\mathcal{T}$ are exactly the same for $\pi(x;q,\mathbf{a})$ and $\pi(x;q,\mathbf{b})$.  Lastly, in arriving at \eqref{2.9} we have paid no attention to error terms, and moreover have used a uniform version of the Hardy-Littlewood conjecture, both in terms of the size of the parameters in the set ${\mathcal A} \cup {\mathcal T}$ (this is relatively minor) and in terms of the size of the set ${\mathcal A} \cup {\mathcal T}$.   To mitigate the last point, we note that in expanding out the inclusion-exclusion product in \eqref{2.5} we may obtain upper and lower 
 bounds by stopping after an odd or an even number of steps (as in Brun's sieve for example); in this manner 
 only a mildly uniform version of the Hardy-Littlewood conjectures seems needed.  For the present we ignore these details, but it would be desirable to place the conjecture \eqref{2.9} on a firmer footing and we intend to return to this in future work.  
 
 With conjecture \eqref{2.9} in hand, we have a conjecture for $\pi(x;q,(a,b))$: namely, we sum the quantity in 
 \eqref{2.9} over all positive integers $h \equiv b-a \pmod q$.  Thus, we expect that 
 \begin{equation} 
 \label{2.10} 
 \pi(x;q,(a,b)) \sim \frac{1}{q} \int_2^x \alpha(y)^{\epsilon_q(a,b)} \Big( \frac{q}{\phi(q) \log y} \Big)^2 {\mathcal D}(a,b;y) dy, 
 \end{equation} 
 say,  where 
 \begin{equation} 
 \label{2.11} 
 {\mathcal D}(a,b;y) =\sum_{\substack{ h> 0\\ h\equiv b-a \pmod q}}  \sum_{{\mathcal A} \subset \{0, h\}} \sum_{\substack{{\mathcal T} \subset [1, h-1]  \\ (t+a,q)=1 \forall t\in{\mathcal T}}} (-1)^{|{\mathcal T}|} {\mathfrak S}_{q,0}({\mathcal A}\cup {\mathcal T}) \Big( \frac{q}{\phi(q)\alpha(y)\log y}\Big)^{|{\mathcal T}|} \alpha(y)^{h\phi(q)/q}.   
 \end{equation} 
 
\noindent {\bf Discarding singular series involving sets with three or more elements.}  We now conjecture that only terms with ${\mathcal A} = {\mathcal T} = \emptyset$ (which gives rise to the main term of $\text{li}(x)/\phi(q)^2$ for $\pi(x;q,(a,b))$), and $|{\mathcal A}| + |{\mathcal T} | =2$ give significant contributions 
 leading to the Main Conjecture, and that all other terms contribute to $\pi(x;q,(a,b))$ an amount $O(x(\log \log x)^2/(\log x)^{3})$.  To argue this, we will use as a guide the work of Montgomery and Soundararajan \eqref{eqn:montgomery-sound} which shows that sums over singular series exhibit square-root cancelation in 
 each variable.  
 
 Suppose for example that ${\mathcal A}=\emptyset$ and $|{\mathcal T}| =\ell \ge 4$ in \eqref{2.11}.  After summing over the variable $h$, these terms may be thought of as 
 $(\log y)^{1-\ell}$ times an average of ${\mathfrak S}_{q,0}({\mathcal T})$ over $\ell$ element sets ${\mathcal T}$ whose elements are all of size about $\log y$.  
 The estimate \eqref{eqn:montgomery-sound} now suggests that this contribution is $\ll (\log \log y)^{\ell/2} (\log y)^{1-\ell/2}$, and since $\ell \ge 4$ the final contribution to 
 $\pi(x;q,(a,b))$ is $O(x(\log \log x)^2/(\log x)^3)$.    If $\ell=3$ then the same argument -- drawing on \eqref{eqn:montgomery-sound} with $k=3$ there, so that the main term there 
 vanishes and the bound is $O(h^{3/2-1/21 +\epsilon})$ -- indicates that such terms contribute to $\pi(x;q,(a,b))$ an amount $O(x (\log x)^{-5/2-1/21+\epsilon})$ which is already smaller than the secondary main terms claimed in the Main Conjecture.  We believe that when $k$ is odd, the work of Montgomery and Soundararajan can be refined and the actual size of the sum in \eqref{eqn:montgomery-sound} is $h^{(k-1)/2}(\log h)^{(k+1)/2}$.  We will pursue this in future work, noting for the present that this expectation suggests that the terms with ${\mathcal A}=\emptyset$ and $|{\mathcal T}|=3$ also make a contribution of $O(x(\log \log x)^2/(\log x)^3)$.
 
 When ${\mathcal A} = \{0\}$ or $\{h\}$, then a similar heuristic to the above shows that terms with $|{\mathcal T}| \ge 2$ make a contribution to $\pi(x;q,(a,b))$ of $O(x(\log \log x)^2/(\log x)^3)$.  Finally if ${\mathcal A}= \{ 0, h\}$ and $|{\mathcal T}| = \ell \ge 1$, then the contribution to \eqref{2.11} may be roughly thought of as $(\log y)^{-\ell}$ times an average of 
 singular series ${\mathfrak S}_{q,0}(\{0\} \cup {\mathcal T}^{+})$ where ${\mathcal T}^+$ (standing for ${\mathcal T} \cup \{h\}$) runs over $\ell +1$ element sets with elements of 
 size $\log y$.  Since the singular series ${\mathfrak S}_{q,0}$ is translation invariant, one can think of this last sum as being $1/(\log y)$ times the average over $\ell+2$ element sets with all elements of size $\log y$.   After making this observation, we can draw on \eqref{eqn:montgomery-sound} (with its proposed refinement for odd $k$) as earlier and this leads to the prediction that the contribution to $\pi(x;q,(a,b))$ of terms with ${\mathcal A}=\{0,h\}$ and any non-empty ${\mathcal T}$ is $O(x(\log \log x)^2/(\log x)^3)$.

 Thus, discarding all terms with $|{\mathcal A}|+ |{\mathcal T}| \ge 3$, we now replace the density ${\mathcal D}(a,b;y)$ in \eqref{2.11} with 
 \begin{equation} 
 \label{2.12} 
 {\mathcal D}(a,b;y)  = {\mathcal D}_0(a,b;y) + {\mathcal D}_1(a,b;y) + {\mathcal D}_2(a,b;y), 
 \end{equation} 
 where (keeping in mind that ${\mathfrak S}_{q,0}$ is $1$ for the empty set, and $0$ for a singleton)
 \begin{equation}
 \label{2.13}  
 \mathcal{D}_0(a,b;y) = \sum_{\substack{ h> 0\\ h\equiv b-a \pmod q}}  (1+ {\mathfrak S}_{q,0}(\{0,h\}) ) \alpha(y)^{h\phi(q)/q}, 
 \end{equation} 
\begin{equation} 
\label{2.14} 
\mathcal{D}_1(a,b;y)=  - \frac{q}{\phi(q)\alpha(y)\log y}  \sum_{\substack{ h> 0\\ h\equiv b-a \pmod q}} 
\sum_{\substack{ t\in [1,h-1] \\ (t+a,q)=1} } ( {\mathfrak S}_{q,0}(\{0,t\}) + {\mathfrak S}_{q,0}(\{ t, h\}) \alpha(y)^{ h\phi(q)/q},  
\end{equation} 
and 
\begin{equation} 
\label{2.15} 
\mathcal{D}_2(a,b;y)= \Big(\frac{q}{\phi(q)\alpha(y)\log y} \Big)^2  \sum_{\substack{ h> 0\\ h\equiv b-a \pmod q}} 
\sum_{\substack{ 1\le t_1 < t_2 < h \\ (t_1+a,q)=(t_2+a,q)=1} } {\mathfrak S}_{q,0}(\{t_1,t_2\}) \alpha(y)^{h\phi(q)/q}.
\end{equation}  
Inserting this in \eqref{2.10}, we thus conjecture that up to $O(x(\log \log x)^2/(\log x)^3)$, there holds 
\begin{equation} 
\label{maineqn} 
\pi(x;q,(a,b)) = \frac{q}{\phi(q)^2} \int_2^x \frac{\alpha(y)^{\epsilon_q(a,b)}}{(\log y)^2} \big( {\mathcal D}_0 +{\mathcal D}_1+{\mathcal D}_2\big)(a,b;y) dy. 
\end{equation} 
 
\noindent{\bf The main proposition.}  To evaluate the sums over two-term singular series above, we invoke the following proposition 
 whose proof we defer to the next section. 
 
\begin{proposition}
\label{prop:sing_sum}  Let $q\geq 2$, and let $v\pmod q$ be any residue class.  
For any positive real number $H$ define 
$$ 
S_0(q,v;H) = \sum_{\substack{h >0 \\ h\equiv v\pmod q} } {\mathfrak S}_{q,0}(\{0, h\}) e^{-h/H}.  
$$   
Then we may write 
 \[
S_0(q,0;H) = - \frac{\phi(q)}{2q} \log H + S_0^c(q,0) + Z_{q,0}(H) + O(H^{-1+\epsilon}), 
\] 
where 
$$ 
S_0^{c}(q,0) = 
 \frac{\phi(q)}{2q}\log \frac{q}{2\pi} -\frac{\phi(q)}{2q}  \sum_{p\mid q} \frac{\log p}{p-1}  +\frac{1}{2},
$$ 
and for any $v\pmod q$, the quantity $Z_{q,v}(H)$ is described in \eqref{eqn:zeros} below, and satisfies the bound $Z_{q,v}(H)=O(H^{-1/2+\epsilon})$, and which we conjecture to be $O(H^{-3/4})$.  
Further, if $(v,q)=d$ with $d<q$, then 
$$ 
S_0(q,v;H) = S_0^c(q,v) + Z_{q,v}(H) + O(H^{-1+\epsilon}), 
$$ 
where 
\begin{align*}
S_0^c(q,v)
	= - \frac{\phi(q)}{2q}\cdot\frac{\Lambda(q/d)}{\phi(q/d)} - B_q(v)
	+ \frac{1}{\phi(q/d)} \sum_{\chi\neq \chi_0\pmod{q/d}} \bar\chi(v/d) L(0,\chi)L(1,\chi) A_{q,\chi},
\end{align*}
with $B_q(v) =\frac 12- \frac{v}{q}$ for $1\le v\le q$ and extended periodically for all $v$, and 
$$
A_{q,\chi} = \prod_{p\mid q} \Big(1-\frac{\chi(p)}{p}\Big) \prod_{p\nmid q} \Big(1-\frac{(1-\chi(p))^2}{(p-1)^2}\Big).
$$
\end{proposition}
 
\noindent{\bf Completing the heuristic.}  Returning to our heuristic calculation, we will apply Proposition \ref{prop:sing_sum} with 
\begin{equation}
\label{eqn:H-def}
H = H(y) :=  -\frac{q}{\phi(q)}\cdot\frac{1}{\log \alpha(y)} = \log y - \frac{q}{2\phi(q)} + O\Big(\frac{1}{\log y}\Big).
\end{equation}
We begin by simplifying a bit the expressions for ${\mathcal D}_0$, ${\mathcal D}_1$ and ${\mathcal D}_2$, discarding terms of size $O(\log \log y/\log y)$ which are negligible for the Main Conjecture.  
 Thus, after summing the geometric series and using \eqref{eqn:H-def},  
\begin{align} 
\label{2.17}
{\mathcal D}_0
&= S_0(q,b-a;H) + \sum_{h\equiv b-a\pmod{q}} e^{-h/H} = S_0(q,b-a;H) + \frac{H}{q} + B_q(b-a) + O\Big(\frac{1}{H}\Big) 
\nonumber \\ 
& = \frac{\log y}{q} + S_0(q,b-a;H) + B_q(b-a) - \frac{1}{2\phi(q)} + O\Big(\frac{1}{\log y}\Big).
\end{align} 

 The definition of ${\mathcal D}_1$ involves two singular series ${\mathfrak S}_{q,0}(\{0,t\})$ and ${\mathfrak S}_{q,0}(t,h)$.  Consider the terms arising from the second case.  Replace $\mathfrak{S}_{q,0}(\{t , h\} )$ by 
 $\mathfrak{S}_{q,0}(\{ 0, r\})$ where $r= h-t$ also lies in $[1,h-1]$ and note that the condition $(t+a,q)=1$ becomes 
 $(r-b,q)=1$.  Thus, ignoring terms of size $O(\log \log y/\log y)$, the second case in ${\mathcal D}_1$ contributes 
 $$ 
 -\frac{q}{\phi(q) \alpha(y) \log y}  \sum_{\substack{r>0 \\ (r-b,q)=1}  } {\mathfrak S}_{q,0} ( \{0, r\}) 
 \sum_{\substack{h>r \\ h\equiv b-a \pmod q}} e^{-h/H} = - \frac{1}{\phi(q)} \sum_{\substack{v\pmod q\\ (v-b,q)=1}} 
S_0(q,v;H).  
 $$ 
 Arguing similarly with the first case, we conclude that 
 \begin{equation} 
 \label{2.18} 
 {\mathcal D}_1 = -\frac{1}{\phi(q)} \sum_{\substack{v\pmod q \\ (v+a,q)=1}} S_0(q,v;H) - 
 \frac{1}{\phi(q)} \sum_{\substack{v\pmod q \\ (v-b,q)=1} } S_0(q,v;H) + O\Big(\frac{\log \log y}{\log y}\Big). 
 \end{equation}

Finally, note that 
\begin{align*}
\sum_{h \equiv b-a\pmod{q}} e^{-h/H} \sum_{\begin{subarray}{c} 1\leq t_1 < t_2 < h \\ (t_1+a,q)=1 \\ (t_2+a,q)=1 \end{subarray}}& \mathfrak{S}_{q,0}(\{t_1,t_2\}) 
	= \sum_{\begin{subarray}{c} 1\leq t_1 < t_2 < h \\ (t_1+a,q)=1 \\ (t_2+a,q)=1 \end{subarray}} \mathfrak{S}_{q,0}(\{0,t_2-t_1\}) \sum_{\begin{subarray}{c} h \equiv b-a\pmod{q} \\ h>t_2 \end{subarray}} e^{-h/H} \\
	&= \frac{H^2}{q^2} \sum_{\begin{subarray}{c} v_1,v_2 \pmod{q} \\ (v_1,q)=1 \\ (v_2,q)=1 \end{subarray}} S_0(q,v_2-v_1;H) + O(H\log H),
\end{align*} 
so that
\begin{equation} 
\label{2.19} 
{\mathcal D}_2 = 
 \frac{1}{\phi(q)^2} \sum_{\begin{subarray}{c} v_1, v_2 \pmod{q} \\ (v_1,q)=1 \\ (v_2,q)=1 \end{subarray}} S_0(q,v_2-v_1;H) + O\Big(\frac{\log\log y}{\log y}\Big).
\end{equation}

Using Proposition \ref{prop:sing_sum} to evaluate \eqref{2.17}, \eqref{2.18}, and \eqref{2.19}, and then inserting 
that in \eqref{2.10} leads to the Main Conjecture.   The term involving $c_1(q;(a,b))$ arises from from terms involving $S_0(q,0;H)$ which has a leading term of size $\log H$ while all other $S_0(q,v;H)$ are only of constant size.  Thus isolating the $-\frac{\phi(q)}{2q} \log H$ leading contribution to $S_0(q,0;H)$ and tracking its appearance in our expressions for ${\mathcal D}_0$, ${\mathcal D}_1$ and ${\mathcal D}_2$ gives 
\begin{align*}
&-\frac{\phi(q)}{2q} (\log H) \delta(a=b) - \frac{2}{\phi(q)} \Big(-\frac{\phi(q)}{2q}  \log H\Big) + \frac{1}{\phi(q)} 
\Big( -\frac{\phi(q)}{2q} \log H\Big)\\
& = \frac{\phi(q)}{2q} (\log \log y) \Big( \frac{1}{\phi(q)} -\delta(a=b)\Big) + O\Big(\frac{\log \log y}{\log y}\Big). 
\end{align*}

The term involving $c_2(q;(a,b))$ is complicated, but follows straightforwardly from our work above.  Having already treated the term $-\frac{\phi(q)}{2q} \log H$ term arising in $S_0(q,0)$, the contributions leading to $c_2(q;(a,b))$ come from the $S_0^c(q,v)$ terms in Proposition \ref{prop:sing_sum}.  
We thus have
\begin{align}
\label{2.20}
\frac{c_2(q;\mathbf{a})}{q}
	&= -\frac{\varepsilon_q(a,b)}{\phi(q)}  + S_0^c(q,b-a)+B_q(b-a)-\frac{1}{2\phi(q)} - \frac{1}{\phi(q)} \sum_{\begin{subarray}{c} v\pmod{q}\\ (v+a,q)=1  \end{subarray}}S_0^c(q,v) \nonumber\\
	&\quad  -\frac{1}{\phi(q)}\sum_{\begin{subarray}{c} v\pmod{q} \\ (v-b,q)=1 \end{subarray}}S_0^c(q,v)+\frac{1}{\phi(q)^2} \sum_{\begin{subarray}{c} v_1, v_2 \pmod{q} \\ (v_1,q)=1 \\ (v_2,q)=1 \end{subarray}} S_0^c(q,v_2-v_1).
\end{align} 
With $C_{q,\chi} = L(0,\chi) L(1,\chi) A_{q,\chi}$ (which is zero unless $\chi$ is an odd character), we may also derive the following alternative expression:
\begin{align} 
\label{2.21}
\frac{c_2(q;\mathbf{a})}{q}	&= \frac{\log 2\pi}{2q} + S_0^c(q,b-a) +B_q(b-a) \nonumber \\ 
&- \frac{1}{\phi(q)} \sum_{\substack{d |q \\ d>1}} \frac{1}{\phi(d)} \sum_{\begin{subarray}{c} \chi \pmod{d}\\ \chi(-1)=-1 \end{subarray}} C_{q,\chi} \Big(
\sum_{\begin{subarray}{c} u\pmod{d} \\  (uq/d+a,q)=1 \end{subarray}} + \sum_{\begin{subarray}{c} u\pmod{d} \\ (uq/d-b,q)=1 \end{subarray}}\Big) \bar\chi(u).
\end{align}

If $\chi$ is induced by the primitive character $\chi^*$, then, writing $\chi= \chi_{0,m} \chi^*$ for some $m$ coprime to the conductor of $\chi^*$, we have 
\[
C_{q,\chi} = C_{q,\chi^*} \prod_{p\mid m} (1-\chi^*(p)).
\]
Further, it is helpful to write $q=q_0 2^r$ with $q_0$ odd.  If now $\chi$ is a character to an odd modulus and $q$ is even, then
\[
C_{q,\chi} = \frac{\bar\chi(2)}{2}C_{q_0,\chi}.
\]
Using these facts, it is possible to simplify the formula in \eqref{2.21} further, and obtain  
\begin{align}
\label{eqn:c2_final}
c_2(q;(a,b)) 
	&= \frac{\log 2\pi}{2} + qS_0^c(q,b-a) +qB_q(b-a) \nonumber \\
	& \quad\quad -\frac{q_0}{\phi(q_0)} \sum_{d\mid q_0} \frac{\mu(d)}{\phi(d)} \sum_{\chi \pmod{d}} C_{q_0,\chi} (\bar\chi(b)-\bar\chi(a)).
\end{align}
For example, if $q$ is prime and $a \neq b$ then 
$$ 
c_2(q;(a,b)) = \frac 12 \log \frac{2\pi}{q} + \frac{q}{\phi(q)} \sum_{\chi \neq \chi_0} C_{q,\chi} \Big( \overline{\chi}(b-a) +
\frac{1}{\phi(q)} (\overline{\chi}(b)-\overline{\chi}(a)) \Big). 
$$

This completes our discussion of the Main Conjecture in the case $r=2$, and the other conjectures follow as simple consequences.

\section{Proof of the Proposition} 
\label{sec:prop} 

\noindent The proof follows along standard lines, and the closely related case of evaluating asymptotically 
$\sum_{h\le H} {\mathfrak S}_0(\{ 0,h\}) (H-h)$ is mentioned in \cite{Goldston} and treated in detail 
in \cite{MS2}.   We will therefore be brief.  Let $\chi$ be a Dirichlet character modulo $m|q$;  possibly $\chi$ could be 
imprimitive, or the principal character.   Define, for Re$(s)>1$, 
\begin{align*}
F_{q,\chi}(s) &:= \sum_{h \geq 1} \frac{\chi(h)}{h^s} \mathfrak{S}_q(\{0,h\}) \\
& = \prod_{p \mid q} \Big(1-\frac{\chi(p)}{p^s}\Big)^{-1} \prod_{p\nmid q}\Big(1-\frac{1}{(p-1)^2}+\frac{\chi(p)}{p^s} \Big(1-\frac{1}{p}\Big)^{-1} 
\Big(1-\frac{\chi(p)}{p^s}\Big)^{-1} \Big),
\end{align*}
so that 
\begin{equation}
\label{3.1}
\sum_{h \geq 1} \chi(h) \mathfrak{S}_q(\{0,h\}) e^{-h/H} = \frac{1}{2\pi i} \int_{(2)} F_{q,\chi}(s) H^s \Gamma(s)\,ds.
\end{equation}
We now note that
\begin{align*}
F_{q,\chi}(s) 
	&= L(s,\chi) \prod_{p\nmid q} \Big(1 - \frac{1}{(p-1)^2}+\frac{\chi(p)}{p^{s-1}(p-1)^2}\Big) \\
	&= L(s,\chi) L(s+1,\chi) \prod_{p\mid q}\Big(1-\frac{\chi(p)}{p^{s+1}}\Big) \prod_{p\nmid q} \Big(1-\frac{(1-\chi(p)/p^s)^2}{(p-1)^2}\Big),\\
\end{align*}
which furnishes a meromorphic continuation of $F_{q,\chi}(s)$ to Re$(s)>-\frac 12$ with possible poles at $s=0$ or $s=1$ in case $\chi$ is principal.  We may also express the above as 
$$
F_{q,\chi}(s)= \frac{L(s,\chi) L(s+1,\chi)}{L(2s+2,\chi^2)}\prod_{p\mid q}\Big(1+\frac{\chi(p)}{p^{s+1}}\Big)^{-1} \prod_{p\nmid q} \Big(1-\frac{1}{(p-1)^2}+\frac{2p\chi(p)}{(p-1)^2(p^{s+1}+\chi(p))}\Big), 
$$
and now the final product above is analytic in Re$(s)>-1$, but for which the line Re$(s)=-1$ forms a 
natural boundary.  

If $\chi$ is non-principal, then by shifting the line of integration to Re$(s)=-\tfrac 12+\epsilon$ we find that 
the quantity in \eqref{3.1} is $L(0,\chi)L(1,\chi) A_{q,\chi} +O(H^{-\frac 12+\epsilon})$, with the main term coming from the pole of $\Gamma(s)$ at $s=0$.  Moreover, we may even shift the line of integration to Re$(s)=-1+\epsilon$ at the cost of picking up residues from the zeros of $L(2s+2,\chi^2)$.   The contribution from these zeros is 
$$ 
Z_{q,\chi}(H) := \sum_{\substack{ \rho, \ \text{Re}(\rho)>0 \\ L(\rho, \chi^2)=0}} \mathop{\text{ Res}}_{s=\rho/2-1} \Big( F_{q,\chi}(s) H^s \Gamma(s)\Big). 
$$
If we suppose that GRH holds for $L(s,\chi^2)$, that its zeros are simple, and that $|L^{\prime}(\rho,\chi^2)|$ is not too small so that (in view of the 
exponential decay of $\Gamma(s)$) the  sum over residues is absolutely convergent,  then we would expect that $Z_{q,\chi}(H) $ is an oscillating term of size $H^{-\frac 34}$.

If $\chi$ is principal, but $m>1$, then $F_{q,\chi}(s)$ has a pole at $s=1$ with residue $\phi(m)/m$, but  there is 
no pole of $F_{q,\chi}$ at $s=0$ since $L(s,\chi_0) = s \Lambda(m) +O(s^2) $ for $s$ near $0$.  Therefore in this 
situation we find 
\[
\sum_{h\geq 1} \chi_0(h) e^{-h/H} \mathfrak{S}_q(\{0,h\}) = \frac{\phi(m)}{m} H - \frac{\phi(q)}{2q} \Lambda(m) + Z_{q,\chi_0}(H) + O(H^{-1+\epsilon}).
\]

Finally if $m=1$ (and $\chi$ is naturally principal) the corresponding $F_{q,\chi}(s)$ has a simple pole at $s=0$ in addition to the pole at $s=1$.  Thus there is a double pole of the integrand in \eqref{3.1}, and computing residues we obtain that 
\[
\sum_{h\geq 1} e^{-h/H} \mathfrak{S}_q(\{0,h\}) = H - \frac{\phi(q)}{2q}\Big[\log 2\pi H + \sum_{p \mid q} \frac{\log p}{p-1} \Big] + Z_{q,\zeta}(H) + O(H^{-1+\epsilon}).
\]

Since
\[
\sum_{h\equiv v\pmod{q}} e^{-h/H} \mathfrak{S}_q(\{0,h\}) = S_0(q,v;H) + \frac{H}{q}+B_q(v)+O\Big(\frac 1H\Big),
\]
our proposition follows, with
\begin{equation}
\label{eqn:zeros}
Z_{q,v}(H) = \frac{1}{\phi(q/d)} \sum_{\chi \pmod{q/d}} \bar\chi(v/d) Z_{q,\chi}(H/d).
\end{equation}
 
\section{Modifications to the heuristics when $r\ge 3$} 
\label{sec:multi-competitor}

\noindent The ideas leading to the general case of the Main Conjecture are similar to those for $r=2$, and so we just give a brief sketch.  For $r\geq 3$ and $\mathbf{a}=(a_1,\dots,a_r)$, we start by writing $\pi(x;q,\mathbf{a})$ as 
\begin{align*}
\sum_{\begin{subarray}{c} n\leq x\\ n\equiv a_1 \pmod{q} \end{subarray}} \sum_{\begin{subarray}{c} h_1,\dots,h_{r-1} >0 \\ h_i \equiv a_{i+1}-a_i \pmod{q} \end{subarray}} \mathbf{1}_\mathcal{P}(n) \prod_{i=1}^{r-1}\Big[& \mathbf{1}_\mathcal{P}(n+h_1+\dots+h_i) \cdot \\ 
&\cdot \prod_{\begin{subarray}{c} 0<t<h_i \\ (t+a_i,q)=1 \end{subarray}} (1-\mathbf{1}_\mathcal{P}(n+h_1+\dots+h_{i-1}+t)) \Big].
\end{align*}
As before, we expand this out, invoke the Hardy-Littlewood conjectures, and then discard all singular series terms except for the empty set and sets with two elements.   This leads to 
\[
\pi(x;q,\mathbf{a}) = \int_2^x \frac{q^{r-1}}{\phi(q)^r} \Big(1-\frac{q}{\phi(q)\log y}\Big)^{\varepsilon_q(\mathbf{a})} \big(\mathcal{D}_0+\mathcal{D}_1+{\mathcal D}_2\big)(\mathbf{a};y) \frac{dy}{(\log y)^r} + O\Big(\frac{x(\log\log x)^2}{\log^3x}\Big),
\]
where $\varepsilon_q(\mathbf{a})=\varepsilon_q(a_1,a_2)+\dots+\varepsilon_q(a_{r-1},a_r)$ and ${\mathcal D}_0$, ${\mathcal D}_1$, and ${\mathcal D}_2$ 
are certain smooth sums of singular series.  For ${\mathcal D}_0$, we have (with $H=H(y)$ as before)
\[
{\mathcal D}_0 = \sum_{\begin{subarray}{c} h_1,\dots,h_{r-1} >0 \\ h_i \equiv a_{i+1}-a_i \pmod{q} \end{subarray}} e^{-(h_1+\dots+h_{r-1})/H} \Big(1+\sum_{0\leq i < j \leq r-1}\mathfrak{S}_{q,0}(\{0,h_{i+1}+\dots+h_j\}) \Big).
\]
Notice that if $j=i+1$ in the inner summation, the resulting expression is $(H/q)^{r-2}$ times the analogous ${\mathcal D}_0$ term in our calculation for $\pi(x;q,(a_j,a_{j+1}))$.  If $j-i > 1$,  we will need to consider sums of the form
\[
S_0^k(q,v;H) := \sum_{h \equiv v\pmod{q}} h^k e^{-h/H} \mathfrak{S}_{q,0}(\{0,h\}),
\]
where $k=j-i-1$.  This can be understood via contour integration as in Proposition \ref{prop:sing_sum}; a key difference is that for $k\ge 1$,  
we have $S_0^k(q,v;H) = O(H^{k-1/2})$ unless $v=0$, in which case $S_0^k(q,0;H) = -\frac{\phi(q)}{2q}\Gamma(k)H^k + O(H^{k-1/2})$.  
Using this to evaluate ${\mathcal D}_0$, we find that it is (up to $O(H^{r-3})$)
\begin{align*}
&\frac{H^{r-1}}{q^{r-1}}+\frac{H^{r-2}}{q^{r-2}}\sum_{i=1}^{r-1} \Big[S_0(q,a_{i+1}-a_i;H)+B_q(a_{i+1}-a_i) + \sum_{k=1}^{r-i-1} \frac{S_0^k(q,a_{i+k+1}-a_i;H)}{k! H^k} \Big] \\
&\sim \frac{H^{r-1}}{q^{r-1}}+\frac{H^{r-2}}{q^{r-2}}\sum_{i=1}^{r-1} \Big[S_0(q,a_{i+1}-a_i;H)+B_q(a_{i+1}-a_i) - \frac{\phi(q)}{2q} \sum_{k=1}^{r-i-1} \frac{\delta(a_i=a_{i+k+1})}{k} \Big],
\end{align*}
and it is this last term which creates the additional bias (in $c_2(q,{\mathbf a})$) against patterns with a non-immediate repetition.

For ${\mathcal D}_1$, up to $O(H^{r-2})$, we obtain a contribution of $(H/q)^{r-1} (1-\frac{\phi(q)}{q}\log y)^{-1}$ times
\begin{align*}
&\sum_{j=1}^{r-1}\Big[\Big(\sum_{(v+a_{j},q)=1} + \sum_{(v-a_{j+1},q)=1} \Big)S_0(q,v;H) + \sum_{k=1}^{j-1}\sum_{(v,q)=1}\frac{S_0^k(q,v-a_{j-k};H)}{k! H^k} \\ 
&\hskip 2 in  + \sum_{k=1}^{r-1-j}\sum_{(v,q)=1}\frac{S_0^k(q,v+a_{j+1+k};H)}{k!H^k} \Big] \\
\sim &\sum_{j=1}^{r-1}\Big(\sum_{(v+a_{j},q)=1}  + \sum_{(v-a_{j+1},q)=1} \Big)S_0(q,v;H) - \frac{\phi(q)}{q}\sum_{k=1}^{r-2}\frac{r-1-k}{k}.
\end{align*}
Finally from ${\mathcal D}_2$ we obtain $(H/q)^r (1-\frac{\phi(q)}{q}\log y)^{-2}$ times
\begin{align*}
\sum_{j=1}^{r-1} \Big(\sum_{\begin{subarray}{c} (v_1,q)=1 \\ (v_2,q)=1 \end{subarray}} S_0(q,v_2-v_1;H) + \Big.&\Big.\sum_{k=1}^{r-1-j}\sum_{\begin{subarray}{c} (v_1,q)=1 \\ (v_2,q)=1 \end{subarray}} \frac{S_0^k(q,v_1+v_2;H)}{k!H^k} \Big) \\
	&\sim (r-1)\sum_{\begin{subarray}{c} (v_1,q)=1 \\ (v_2,q)=1 \end{subarray}} S_0(q,v_2-v_1;H) - \frac{\phi(q)^2}{2q}\sum_{k=1}^{r-2}\frac{r-1-k}{k}.
\end{align*}
Assembling these contributions yields the Main Conjecture.

\section{Comparison of the conjecture with numerical data} 

\noindent We begin by comparing the Main Conjecture against the data for $r=2$ and $q=3$ or $4$.  In each of these cases, our conjecture is that
\begin{equation}
\label{eqn:q=3or4}
\pi(x;q,\mathbf{a}) = \frac{\mathrm{li}(x)}{4} \Big(1 \pm \frac{1}{2\log x} \log\Big(\frac{2\pi \log x}{q}\Big)\Big) + O\Big(\frac{x}{(\log  x)^{11/4}}\Big),
\end{equation}
with the sign being negative if $a_1 \equiv a_2 \pmod{q}$ and positive if not.  However, in order to obtain \eqref{eqn:q=3or4} in such a clean form, a number of asymptotic approximations were used throughout Section \ref{sec:heuristic}, and it is reasonable to expect that the unsimplified integral expression \eqref{maineqn} for $\pi(x;q,\mathbf{a})$ would provide a better fit to the data.  Indeed, we find the following.

{\small
\begin{center}
\begin{tabular}{l}
\\
Actual \\
\eqref{maineqn} \\
\eqref{eqn:q=3or4} \\
\\
\\
\\
\\
\\
\\
\\
\\
\\
\end{tabular}
\begin{tabular}{|l|l|l|}
\hline $x$ & $\pi(x;3,(1,1))$ & $\pi(x;3,(1,2))$  \\ \hline
$10^9$ 	
	& $1.132 \cdot 10^7$ & $1.411 \cdot 10^7$ \\ 
	& $1.137 \cdot 10^7$ & $1.405 \cdot 10^7$ \\
 	& $1.156 \cdot 10^7$ & $1.387 \cdot 10^7$ \\ \hline
$10^{10}$ 
	& $1.024 \cdot 10^8$ & $1.251 \cdot 10^8$ \\ 
	& $1.028 \cdot 10^8$ & $1.247 \cdot 10^8$ \\ 
 	& $1.042 \cdot 10^8$ & $1.233 \cdot 10^8$ \\ \hline
$10^{11}$ 
	& $9.347\cdot 10^8$ & $1.124\cdot 10^9$ \\ 
	& $9.383\cdot 10^8$ & $1.121\cdot 10^9$ \\
 	& $9.488\cdot 10^8$ & $1.110\cdot 10^9$ \\ \hline
 $10^{12}$ 
 	& $8.600\cdot 10^9$ & $1.020\cdot 10^{10}$ \\ 
 	& $8.630\cdot 10^9$ & $1.017\cdot 10^{10}$ \\
 	& $8.712\cdot 10^9$ & $1.009\cdot 10^{10}$ \\ \hline
\end{tabular} \,
\begin{tabular}{|l|l|}
\hline $\pi(x;4,(1,1))$ & $\pi(x;4,(1,3))$ \\ \hline
$1.141 \cdot 10^7$ & $1.401 \cdot 10^7$ \\ 
$1.148 \cdot 10^7$ & $1.395 \cdot 10^7$ \\
$1.164 \cdot 10^7$ & $1.378 \cdot 10^7$ \\ \hline
$1.032 \cdot 10^8$ & $1.244 \cdot 10^8$ \\ 
$1.037 \cdot 10^8$ & $1.239 \cdot 10^8$ \\ 
$1.049 \cdot 10^8$ & $1.226 \cdot 10^8$ \\ \hline
$9.412 \cdot 10^8$ & $1.118 \cdot 10^9$ \\ 
$9.450 \cdot 10^8$ & $1.114 \cdot 10^9$ \\
$9.547 \cdot 10^8$ & $1.104 \cdot 10^9$ \\ \hline
$8.654 \cdot 10^9$ & $1.015 \cdot 10^{10}$ \\ 
$8.684 \cdot 10^9$ & $1.012 \cdot 10^{10}$ \\
$8.760 \cdot 10^9$ & $1.004 \cdot 10^{10}$ \\ \hline
\end{tabular}
\end{center}
}
\smallskip

Going forward, we will present only the comparison of $\pi(x;q,\mathbf{a})$ against \eqref{maineqn}, so we explain briefly how we compute this approximation.  
In \eqref{2.17}, \eqref{2.18} and \eqref{2.19}, we determined ${\mathcal D}_0$, ${\mathcal D}_1$ and ${\mathcal D}_2$ in terms of $S_0(q,v;H)$ and in the process 
replaced geometric progressions in $h$ with suitable approximations.  Of course the geometric progressions could just be computed exactly.  We keep the exact but messy 
expressions so obtained, and for $S_0(q,v;H)$ use the main terms described in Proposition \ref{prop:sing_sum}.   This yields an expression for $\pi(x;q,\mathbf{a})$ as an explicit integral, which we computed numerically in Sage.  The actual values of $\pi(x;q,\mathbf{a})$ were computed in \verb^C++^ using the primesieve library.  Code for both computations can be found on the first author's website.

Next we consider $q=8$.  Here too the constants simplify, with $c_2(8;(a,b))$ depending only on the difference $b-a \pmod{8}$ (a fact reflected in the data, as predicted by Conjecture \ref{conj:primepower}).  Explicitly, we have $c_2(8;(a,a)) = (5\log 2 - 3\log \pi)/2$, $c_2(8;(a,a+2))=c_2(8;(a,a+6))= (\log \pi - \log 2)/2$, and $c_2(8;(a,a+4)) = (\log \pi - 3 \log 2)/2$.  Thus, we should expect that, among the non-diagonal patterns, those with $b-a=4$ should be the least frequent, and that those with $b-a=2$ and $6$ should be rather close.  Indeed, we find:

{\small
\begin{center}
\begin{tabular}{l}
\\
Actual \\
\eqref{maineqn} \\
\\
\\
\\
\\
\\
\\
\end{tabular}
\begin{tabular}{|l|l|l|l|l|}
\hline $x$ & $\pi(x;8,(1,1))$ & $\pi(x;8,(1,3))$ & $\pi(x;8,(1,5))$ & $\pi(x;8,(1,7))$ \\ \hline
$10^9$ 	
	& $2.356 \cdot 10^6$ & $3.496 \cdot 10^6$ & $3.351 \cdot 10^6$ & $3.508 \cdot 10^6$ \\ 
 	& $2.369 \cdot 10^6$ & $3.462 \cdot 10^6$ & $3.370 \cdot 10^6$ & $3.511 \cdot 10^6$ \\ \hline
$10^{10}$ 
	& $2.170 \cdot 10^7$ & $3.101 \cdot 10^7$ & $2.988 \cdot 10^7$ & $3.117 \cdot 10^7$ \\ 
 	& $2.179 \cdot 10^7$ & $3.081 \cdot 10^7$ & $3.004 \cdot 10^7$ & $3.112 \cdot 10^7$ \\ \hline
$10^{11}$ 
	& $2.010 \cdot 10^8$ & $2.787 \cdot 10^8$ & $2.696 \cdot 10^8$ & $2.802 \cdot 10^8$ \\ 
 	& $2.016 \cdot 10^8$ & $2.775 \cdot 10^8$ & $2.709 \cdot 10^8$ & $2.795 \cdot 10^8$ \\ \hline
 $10^{12}$ 
 	& $1.871 \cdot 10^9$ & $2.530 \cdot 10^9$ & $2.456 \cdot 10^9$ & $2.545 \cdot 10^9$ \\ 
 	& $1.876 \cdot 10^9$ & $2.523 \cdot 10^9$ & $2.466 \cdot 10^9$ & $2.537 \cdot 10^9$ \\ \hline
\end{tabular} 
\end{center}
}
\smallskip

We now turn to the patterns $\pmod{12}$.  Here, the quadratic character $\chi \pmod{3}$ plays a role for those patterns $(a,b)$ with $a\not\equiv b\pmod{3}$.  In particular, it does not play a role in the diagonal patterns, for which $c_2(12;\mathbf{a})$ is given by \eqref{eqn:diagonal_constant}.  For non-diagonal patterns, we have:

\begin{center}
\begin{tabular}{|l|l|l|l|l|}
\hline
$\mathbf{a}$ & $(1,5)$ & $(1,7)$ & $(1,11)$ & $(5,1)$ \\ \hline
$c_2(12;\mathbf{a})$ 
	& $\frac{1}{2}\log(2\pi/9) + \frac{\pi}{\sqrt{3}}A_{12,\chi}$ 
	& $\frac{1}{2}\log(\pi/8)$
	& $ \frac{1}{2}\log(2\pi) - \frac{\pi}{\sqrt{3}}A_{12,\chi}$ 
	& $\frac{1}{2}\log(2\pi/9) - \frac{\pi}{\sqrt{3}}A_{12,\chi}$ \\ \hline
\hline
$\mathbf{a}$ & $(5,7)$ & $(7,1)$ & $(7,5)$ & $(11,1)$ \\ \hline
$c_2(12;\mathbf{a})$ 
	& $\frac{1}{2}\log(2\pi) + \frac{\pi}{\sqrt{3}}A_{12,\chi}$ 
	& $\frac{1}{2}\log(\pi/8)$
	& $ \frac{1}{2}\log(2\pi) - \frac{\pi}{\sqrt{3}}A_{12,\chi}$ 
	& $\frac{1}{2}\log(2\pi) + \frac{\pi}{\sqrt{3}}A_{12,\chi}$ \\ \hline
\end{tabular}
\smallskip 

(The other values of $c_2(12;\mathbf{a})$ are determined by $c_2(12;\mathbf{a}^\mathrm{opp})$.)
\end{center}
\bigskip

 Here, $A_{12,\chi} \approx 1.036$, so that $c_2(12;(5,7))$ and $c_2(12;(11,1))$ are the largest of these.  Moreover, as in the $\pmod{8}$ case, there are symmetries between patterns with the same difference $b-a$.  We find the following.
\bigskip

{\small
\begin{center}
\begin{tabular}{l}
\\
Actual \\
\eqref{maineqn} \\
\\
\\
\end{tabular}
\begin{tabular}{|l|l|l|l|l|l|}
\hline $x$ & $\pi(x;12,(1,1))$ & $\pi(x;12,(1,5))$ & $\pi(x;12,(1,7))$ & $\pi(x;12,(1,11))$ & $\pi(x;12,(5,1))$ \\ \hline
$10^9$ 	
	& $2.305 \cdot 10^6$ & $3.809 \cdot 10^6$ & $3.352 \cdot 10^6$ & $3.245 \cdot 10^6$ & $2.994 \cdot 10^6$ \\ 
 	& $2.364 \cdot 10^6$ & $3.682 \cdot 10^6$ & $3.318 \cdot 10^6$ & $3.347 \cdot 10^6$ & $3.073 \cdot 10^6$ \\ \hline
 $10^{12}$ 
 	& $1.842 \cdot 10^9$ & $2.670 \cdot 10^9$ & $2.458 \cdot 10^9$ & $2.402 \cdot 10^9$ & $2.271 \cdot 10^9$ \\ 
 	& $1.863 \cdot 10^9$ & $2.651 \cdot 10^9$ & $2.448 \cdot 10^9$ & $2.440 \cdot 10^9$ & $2.307 \cdot 10^9$ \\ \hline 	
\end{tabular} 
\bigskip

\hspace{0.52in}
\begin{tabular}{|l|l|l|l|l|l|}
\hline $x$ & $\pi(x;12,(5,5))$ & $\pi(x;12,(5,7))$ & $\pi(x;12,(7,1))$ & $\pi(x;12,(7,5))$ & $\pi(x;12,(11,1))$ \\ \hline
$10^9$ 	
	& $2.305 \cdot 10^6$ & $4.061 \cdot 10^6$ & $3.351 \cdot 10^6$ & $3.245 \cdot 10^6$ & $4.061 \cdot 10^6$ \\ 
 	& $2.365 \cdot 10^6$ & $3.956 \cdot 10^6$ & $3.318 \cdot 10^6$ & $3.347 \cdot 10^6$ & $3.956 \cdot 10^6$ \\ \hline
 $10^{12}$ 
 	& $1.842 \cdot 10^9$ & $2.831 \cdot 10^9$ & $2.458 \cdot 10^9$ & $2.402 \cdot 10^9$ & $2.831 \cdot 10^9$ \\ 
 	& $1.862 \cdot 10^9$ & $2.784 \cdot 10^9$ & $2.448 \cdot 10^9$ & $2.440 \cdot 10^9$ & $2.784 \cdot 10^9$ \\ \hline 	
\end{tabular} 
\end{center}
}
\bigskip

We close by considering $q=5$ (which amounts to considering the last decimal digit of primes).  Essentially no simplfications can be made for the constants $c_2(q;\mathbf{a})$.  For any non-diagonal pattern $(a,b)$, we find
\begin{equation*}
\label{eqn:c2-5}
c_2(5;(a,b)) = \frac{\log (2\pi/5)}{2} + \frac{5}{2} \text{Re }\Big(L(0,\chi)L(1,\chi)A_{5,\chi}\Big[\bar\chi(b-a)+\frac{\bar\chi(b)-\bar\chi(a)}{4}\Big]\Big),
\end{equation*}
where $\chi$ is either of the complex characters $\pmod{5}$.  Apart from the understood symmetry $c_2(5;(a,b))=c_2(5;(-b,-a))$, the value of $c_2$ determines the pattern.  Thus, we might expect significant variation between the various patterns, and in particular no additional symmetries like we saw $\pmod{8}$ and $\pmod{12}$.  We find, presenting only the 
first of $(a,b)$ and $(-b,-a)$:

{\small
\begin{center}
\begin{tabular}{l}
\\
Actual \\
\eqref{maineqn} \\
\\
\\
\end{tabular}
\begin{tabular}{|l|l|l|l|l|l|}
\hline $x$ & $\pi(x;5,(1,1))$ & $\pi(x;5,(1,2))$ & $\pi(x;5,(1,3))$ & $\pi(x;5,(1,4))$ & $\pi(x;5,(2,1))$ \\ \hline
$10^9$ 	
	& $2.328 \cdot 10^6$ & $3.842 \cdot 10^6$ & $3.796 \cdot 10^6$ & $2.745 \cdot 10^6$ & $3.244 \cdot 10^6$ \\ 
 	& $2.354 \cdot 10^6$ & $3.774 \cdot 10^6$ & $3.835 \cdot 10^6$ & $2.750 \cdot 10^6$ & $3.149 \cdot 10^6$ \\ \hline
 $10^{12}$ 
 	& $1.848 \cdot 10^9$ & $2.704 \cdot 10^9$ & $2.706 \cdot 10^9$ & $2.145 \cdot 10^9$ & $2.386 \cdot 10^9$ \\ 
 	& $1.863 \cdot 10^9$ & $2.682 \cdot 10^9$ & $2.717 \cdot 10^9$ & $2.141 \cdot 10^9$ & $2.352 \cdot 10^9$ \\ \hline 	
\end{tabular} 
\bigskip

\hspace{0.52in}
\begin{tabular}{|l|l|l|l|l|l|}
\hline $x$ & $\pi(x;5,(2,2))$ & $\pi(x;5,(2,3))$ & $\pi(x;5,(3,1))$ & $\pi(x;5,(3,2))$ & $\pi(x;5,(4,1))$ \\ \hline
$10^9$ 	
	& $2.228 \cdot 10^6$ & $3.444 \cdot 10^6$ & $3.047 \cdot 10^6$ & $3.595 \cdot 10^6$ & $4.092 \cdot 10^6$ \\ 
 	& $2.337 \cdot 10^6$ & $3.391 \cdot 10^6$ & $3.033 \cdot 10^6$ & $3.568 \cdot 10^6$ & $4.176 \cdot 10^6$ \\ \hline
 $10^{12}$ 
 	& $1.811 \cdot 10^9$ & $2.499 \cdot 10^9$ & $2.301 \cdot 10^9$ & $2.586 \cdot 10^9$ & $2.867 \cdot 10^9$ \\ 
 	& $1.856 \cdot 10^9$ & $2.477 \cdot 10^9$ & $2.295 \cdot 10^9$ & $2.570 \cdot 10^9$ & $2.893 \cdot 10^9$ \\ \hline 	
\end{tabular} 
\end{center}
}
\smallskip

An interesting feature to be observed here is that, initially, $\pi(x;5,(1,2))$ is larger than $\pi(x;5,(1,3))$, despite our conjecture predicting the opposite ordering.  In fact, this is true for all $x$ between $41{,}231$ and $5.076 \cdot 10^{11}$.  However, at about $5.082 \cdot 10^{11}$, $\pi(x;5,(1,3))$ becomes consistently larger, seemingly forever, exactly as our conjecture would predict.  We take this as reasonable evidence for our speculation that there are even more lower-order terms (e.g., on the order of $x(\log\log x)^2/\log^3x$), which in this case apparently conspire to point in the opposite direction than the bias in the Main Conjecture.

\bibliographystyle{abbrv}
\bibliography{primebias}

\end{document}